\documentclass[10pt]{amsproc} 
\usepackage[all]{xy}
\SelectTips{eu}{}
\usepackage{ifthen} 
\usepackage{amssymb}
\calclayout
\makeatletter 
\def\serieslogo@{} 
\def\@setcopyright{} 
\makeatother
\title[Exercises on derived categories]{Exercises 
on derived categories, resolutions, and Brown representability}
\author[Henning Krause]{Henning Krause}
\address{Henning Krause\\ Institut f\"ur Mathematik\\
Universit\"at Paderborn\\ 33095 Paderborn\\ Germany.}
\email{hkrause@math.upb.de}

\newcommand{\smatrix}[1]{\left[\begin{smallmatrix}#1\end{smallmatrix}\right]}

\newcommand{\Mod}{\operatorname{Mod}\nolimits}

\newcommand{\End}{\operatorname{End}\nolimits}
\newcommand{\Hom}{\operatorname{Hom}\nolimits}
\newcommand{\RHom}{\operatorname{{\bfR}Hom}\nolimits}
\newcommand{\END}{\operatorname{\mathcal E\!\!\:\mathit n\mathit d}\nolimits}
\newcommand{\HOM}{\operatorname{\mathcal H\!\!\:\mathit o\mathit m}\nolimits}

\newcommand{\Ker}{\operatorname{Ker}\nolimits}

\newcommand{\Ext}{\operatorname{Ext}\nolimits}

\newcommand{\Inj}{\operatorname{Inj}\nolimits}
\newcommand{\Proj}{\operatorname{Proj}\nolimits}
 
\newcommand{\op}{\mathrm{op}}
 
\newcommand{\id}{\mathrm{id}}
\newcommand{\dg}{\mathrm{dg}}

\newcommand{\inj}{\mathrm{inj}}

\newcommand{\comp}{\mathop{\raisebox{+.3ex}{\hbox{$\scriptstyle\circ$}}}}
\newcommand{\lto}[1][{}]{\stackrel{#1}{\longrightarrow}} 
\renewcommand{\to}[1][{}]{\stackrel{#1}{\rightarrow}}

\def\p{\phi}

\def\La{\Lambda}
\def\Si{\Sigma}

\def\A{{\mathcal A}}

\def\I{{\mathcal I}}

\def\P{{\mathcal P}}

\def\T{{\mathcal T}}
\def\U{{\mathcal U}}

\def\bbZ{\mathbb Z}

\def\bfC{\mathbf C}
\def\bfD{\mathbf D}

\def\bfK{\mathbf K}

\def\bfp{\mathbf p}

\def\bfR{\mathbf R}
\def\bfS{\mathbf S}

\begin{document}
\maketitle

\noindent The numbering of the following exercises refers to the article
``Derived categories, resolutions, and Brown representability'' in
this volume.

\bigskip
\noindent{\bf (1.2.1)} Let $\A$ be an abelian category.  Show that
$\bfK(\A)$ and $\bfD(\A)$ are additive categories and
that the canonical functor $\bfK(\A)\to\bfD(\A)$ is additive.

\bigskip
\noindent{\bf (1.4.1)} Let $\A$ be an abelian category and denote by
$T$ the class of all quasi-isomorphisms in $\bfC(\A)$.  Show that two
maps $\p,\psi\colon X\to Y$ in $\bfC(\A)$ are identified by the
canonical functor $\bfC(\A)\to\bfC(\A)[T^{-1}]$ if $\p-\psi$ is
null-homotopic.

\bigskip
\noindent{\bf (1.5.1)} Let $\A$ be the module category of a ring
$\La$. Show that $\Hom_{\bfD(\A)}(\La,X)\cong H^0X$ for every
complex $X$ of $\La$-modules.

\bigskip
\noindent{\bf (1.5.2)} Let $\A$ be an abelian category. Show that the
canonical functor $\A\to\bfD(\A)$ identifies $\A$ with the full
subcategory of complexes $X$ in $\bfD(\A)$ such that $H^nX=0$ for all
$n\neq 0$.

\bigskip
\noindent{\bf (1.6.1)} Let $\A$ be the category of vector spaces over a
field $k$. Describe all objects and morphisms in $\bfD(\A)$.

\bigskip
\noindent{\bf (1.6.2)} Let $\A$ be the category of finitely generated
abelian groups and $\P$ be the category of finitely generated free
abelian groups. Describe all objects and morphisms in
$\bfD^b(\A)$. Show that the canonical functor
$\bfK^b(\P)\to \bfD^b(\A)$ is an equivalence.

\bigskip
\noindent{\bf (1.6.3)} Let $k$ be a field and consider the following finite
dimensional algebras.
$$\La_1=\smatrix{k&k&k\\0&k&k\\0&0&k}\quad
\La_2=\smatrix{k&k&0\\0&k&0\\0&k&k}\quad\La_3=\La_1/I,\;\; I=
\smatrix{0&0&k\\0&0&0\\0&0&0}$$ Describe in each case the category
$\A_i$ of finite dimensional $\La_i$-modules and its derived category
$\bfD^b(\A_i)$.  Here are some hints.
\begin{enumerate}
\item $\A_1$ and $\A_2$ are hereditary categories, but $\A_3$ is not.
\item Each object in $\A_i$ or $\bfD^b(\A_i)$ decomposes
essentialy uniquely into a finite number of indecomposable objects.
\item The indecomposable projective $\La_i$-modules are $E_{jj}\La_i$,
$j=1,2,3$.
\item $\La_1$ and $\La_2$ have each $6$ pairwise non-isomorphic
indecomposable modules, and $\La_3$ has $5$.
\item $\Ext^n_{\La_i}(X,Y)$ has $k$-dimension at most $1$ for all
indecomposable $\La_i$-modules $X,Y$ and $n\geq 0$.
\end{enumerate}
The {\em Auslander-Reiten quiver} provides a convenient method to display
the categories $\A_i$ and $\bfD^b(\A_i)$, because the
morphism spaces between indecomposable objects are at most
one-dimensional. This quiver (=oriented graph) is defined as follows.
The vertices correspond to the indecomposable objects. Put an arrow
$X\to Y$ between two indecomposable objects if there is an irreducible
map $\p\colon X\to Y$ (where $\p$ is {\em irreducible} if $\p$ is not
invertible and any factorization $\p=\p''\comp\p'$ implies that $\p'$
is a split monomorphism or $\p''$ is a split epimorphism).

\bigskip
\noindent{\bf (1.7.1)} Let $\A$ be an abelian category. Show that the
canonical functor $\bfD^b(\A)\to\bfD(\A)$ is fully faithful.

\bigskip
\noindent{\bf (1.7.2)} Let $\A$ be an abelian category and denote by
$\I$ the full subcategory of injective objects. Suppose that $\A$ has
enough injective objects. Then the canonical functor
$\bfK^+(\I)\to\bfD^+(\A)$ is an equivalence.

\bigskip
\noindent{\bf (1.7.3)} Let $\A$ be the category of finite dimensional
modules over $\Lambda=k[T]/(T^2)$, where $k$ is a field. Describe the
derived category $\bfD^b(\A)$. (Hint: Fix an injective resulution $I$
of the unique simple module $k[T]/(T)$ (with $I^n=\La$ or $I^n=0$ for
all $n$) and build every object in $\bfD^b(\A)$ from $I$.)

\bigskip
\noindent{\bf (2.1.1)} Let $\T$ be a triangulated category. Show that the
coproduct of two exact triangles is an exact triangle. Generalize this
as follows. Let $X_i\to Y_i\to Z_i\to \Si X_i$ be a family of exact
triangles such that the coproducts $\coprod_iX_i$, $\coprod_iY_i$, and
$\coprod_iX_i$ exist in $\T$. Show that
$$\coprod_iX_i\lto \coprod_iY_i\lto \coprod_iZ_i\lto \Si
\big(\coprod_i X_i\big)$$ is an exact triangle in $\T$.

\bigskip
\noindent{\bf (2.1.2)} Let $\T$ be a triangulated category. Show that
the opposite category $\T^\op$ is also triangulated.

\bigskip
\noindent{\bf (2.3.1)} Show that every monomorphism $\p\colon X\to Y$
in a triangulated category has a left inverse $\p'$ such that
$\p'\comp\p=\id_X$.

\bigskip
\noindent{\bf (2.4.1)} Give an example of an exact triangle
$\Delta$ and two endomorphisms $(\p_1,\p_2,\p_3')$ and
$(\p_1,\p_2,\p_3'')$ of $\Delta$ such that
$\p_3'\neq\p_3''$.

\bigskip
\noindent{\bf (2.5.1)} Let $\A$ be an additive category. Check the
axioms (TR1) -- (TR4) for $\bfK(\A)$.

\bigskip
\noindent{\bf (3.1.1)} Let $\A$ be an abelian category. Show that a
map in $\bfK(\A)$ is a quasi-isomorphism if and only if the canonical
functor $\bfK(\A)\to\bfD(\A)$ sends the map to an isomorphism in
$\bfD(\A)$.

\bigskip
\noindent{\bf (3.2.1)} Let $F\colon\T\to\U$ be an exact functor
between triangulated categories. Show that a right adjoint of $F$ is
an exact functor.

\bigskip
\noindent{\bf (3.2.2)} Let $\A$ be an abelian category. Find a
criterion such that $\bfD(\A)$ is an abelian category.

\bigskip
\noindent{\bf (3.3.1)} Let $\La$ be a noetherian ring and $\A$ be the
category of $\La$-modules. A complex $X$ in $\A$ has {\em finite
cohomology} if $H^nX$ is finitely generated for all $n$ and vanishes
for almost all $n\in\bbZ$ . Show that the complexes with finite
cohomology form a thick subcategory of $\bfD(\A)$.

\bigskip
\noindent{\bf (3.3.2)} Let $\A$ be the category of finite dimensional
modules over $k[T]/(T^n)$. Describe the thick subcategory of all
acyclic complexes in $\bfK(\A)$ which have projective components. Draw
the Auslander-Reiten quiver of this category.  (Hint: Note that
projective and injective modules over $k[T]/(T^n)$ coincide. Each
acyclic complex $X$ of injectives is essentially determined by the
module $Z^0X$.)

\bigskip
\noindent{\bf (3.5.1)} Let $\La$ be a ring and $e=e^2\in\La$ be an
idempotent. Let $\Gamma=e\La e\cong\End_\La(e\La)$. Then
$\Hom_\La(e\La,-)$ induces an exact functor $\Mod\La\to\Mod\Gamma$
which extends to an exact functor $F\colon
\bfD(\Mod\La)\to\bfD(\Mod\Gamma)$. Show that $F$ induces an
equivalence
$$\bfD(\Mod\La)/{\Ker F}\to\bfD(\Mod\Gamma).$$

\bigskip
\noindent{\bf (4.1.1)} Let $\A$ be an additive category. Give a
presentation of the cokernel of a map between two coherent functors in
$\widehat\A$.

\bigskip
\noindent{\bf (4.1.2)} Let $\A$ be an additive category.  Show that
for every family of functors $F_i$ in $\widehat\A$ having a
presentation
$$\A(-,X_i)\lto[(-,\p_i)] \A(-,Y_i)\lto F_i\lto 0,$$
the coproduct $F=\coprod_i F_i$ in $\widehat\A$ has a presentation
$$\A(-,\coprod_iX_i)\lto[(-,\amalg\p_i)] \A(-,\coprod_iY_i)\lto F\lto 0.$$

\bigskip
\noindent{\bf (4.1.3)} Let $\La$ be a ring and $\A$ be the category of
free $\La$-modules. Show that $\widehat\A$ is equivalent to the
category of $\La$-modules.

\bigskip
\noindent{\bf (4.2.1)} Let $F\colon\T\to\U$ be an exact functor
between triangulated categories. Show that the induced functor
$\widehat\T\to\widehat\U$ is exact.

\bigskip
\noindent{\bf (4.5.1)} Let $\A$ be the category of $\La$-modules over
a ring $\La$. Show that $\La$ is a perfect generator for $\bfD(\A)$.

\bigskip
\noindent{\bf (4.5.2)} Let $\T$ be a triangulated category with
arbitrary coproducts. Show that one can replace in the definition of a perfect generator
the condition 
\begin{enumerate}
\item[(PG1)] There is no proper full triangulated subcategory of $\T$ which
contains $S$ and is closed under taking coproducts.
\end{enumerate}
by the following condition
\begin{enumerate}
\item[(PG1')] Let $X$ be in $\T$ and suppose $\Hom_\T(\Si^nS,X)=0$ for
all $n\in\bbZ$. Then $X=0$.
\end{enumerate}

\bigskip
\noindent{\bf (5.1.1)} Let $\A$ be an abelian category and $I$
be the injective resolution of an object $A$. Show that the canonical map
$A\to I$ induces an isomorphism
$$\Hom_{\bfK(\A)}(I,X)\cong\Hom_{\bfK(\A)}(A,X)$$
for every complex $X$ with injective components.

\bigskip
\noindent{\bf (5.1.2)} Let $\A$ be an abelian category and suppose
$\A$ has arbitrary products. Then the canonical functor
$\bfK(\A)\to\bfD(\A)$ preserves products if and only if products in
$\A$ are exact.

\bigskip
\noindent{\bf (5.1.3)} Let $\A$ be an abelian category with a
projective generator. Show that products in $\A$ are exact.

\bigskip
\noindent{\bf (5.1.4)} Let $\A$ be an abelian category with arbitrary
products, and denote by $\Inj\A$ the full subcategory of injective
objects. Show that
$$\bfK^+(\Inj\A)\subseteq\bfK_\inj(\A)\subseteq\bfK(\Inj\A).$$ (Hint: Write
every complex in $\bfK^+(\Inj\A)$ as a homotopy limit of truncations
from $\bfK^b(\Inj\A)$.)

\bigskip
\noindent{\bf (5.1.5)} Let $\A$ be an abelian category with exact
products and an injective cogenerator.  Denote by $\Inj\A$ the full
subcategory of injective objects. Suppose every object in $\A$ has
finite injective dimension.  Show that
$\bfK_\inj(\A)=\bfK(\Inj\A)$. In particular, $\bfK(\Inj\A)$ and
$\bfD(\A)$ are equivalent. (Hint: An acyclic complex of injectives is
null-homotopic.)

\bigskip
\noindent{\bf (5.1.6)} If a ring $\La$ has finite global dimension,
then $\bfK(\Inj\La)$ and $\bfK(\Proj\La)$ are equivalent. 

\bigskip
\noindent{\bf (5.3.1)} Consider the setup from (1.6.3). Define
$\La_1$-modules $$B=E_{11}\La_1\amalg E_{22}\La_1\amalg (E_{22}\La_1/
E_{23}\La_1)\quad\textrm{and}\quad C=(E_{11}\La_1/ E_{12}\La_1)\amalg
E_{11}\La_1\amalg E_{33}\La_1.$$ Show that $\La_2\cong\End_{\La_1}(B)$
and $\La_3\cong\End_{\La_1}(C)$. Viewing these isomorphisms as
identifications, we have bimodules
$_{\La_2}B_{\La_1}$ and $_{\La_3}C_{\La_1}$ which induce equivalences
$$\RHom_{\La_1}(B,-)\colon\bfD^b(\A_1)\to\bfD^b(\A_2)
\quad\textrm{and}\quad
\RHom_{\La_1}(C,-)\colon\bfD^b(\A_1)\to\bfD^b(\A_3).$$ (The
$\La_1$-modules $B$ and $C$ are examples of so-called tilting
modules.)

\bigskip
\noindent{\bf (6.1.1)} Let $k$ be a field and consider again the algebra
$$\La=\smatrix{k&k&k\\0&k&k\\0&0&k}.$$ Denote by $S=S_1\amalg
S_2\amalg S_3$ the coproduct of the three simple $\La$-modules.  Let
$P=\bfp S$ be a projective resolution of $S$. Compute $A=\END_\La(P)$
and show that $H^nA\cong\Ext^n_\La(S,S)$ for all $n$. Show that
$X\mapsto\HOM_\La(P,X)$ induces a functor
$\bfK(\Proj\La)\to\bfD_\dg(A)$ which is an equivalence.

\bigskip
\noindent{\bf (6.2.1)} View a $k$-algebra $A$ as a category $\A$ with
a single object $*$ and $\A(*,*)=A$. Establish an equivalence between
the category of right $A$-modules and the category of $k$-linear
functors $\A^\op\to\Mod k$.

\bigskip
\noindent{\bf (6.5.1)} Let $\A$ be the module category of a noetherian
ring, and let $A$ in $\A$ be finitely generated. Show that $A$ is a
compact object in $\A$. The object $A$ is compact in $\bfD(\A)$ if and
only if $A$ has finite projective dimension.

\bigskip
\noindent{\bf (6.5.2)} Let $\A$ be the module category of a
commutative noetherian ring $\La$. Show that a complex $X$ in
$\bfD(\A)$ has finite cohomology if and only if
$\Hom_{\bfD(\A)}(\Si^nC,X)$ is finitely generated over $\La$ for every
compact object $C$ and all $n\in\bbZ$, and if it vanishes for almost
all $n\in\bbZ$.

\bigskip
\noindent{\bf (7.4.1)} Let $\A$ be an additive category. Show that the
two triangulated structures on $\bfK(\A)$ (defined via mapping cones sequences
and via degree-wise split exact sequences) coincide.

\bigskip
\noindent{\bf (7.4.2)} Let $\La$ be a ring such that projective and
injective $\La$-modules coincide. Then $\La$ is noetherian and the category
$\A$ of finitely generated $\La$-modules is an abelian Frobenius
category. Denote by $\bfD^b(\Proj\A)$ the thick subcategory of
$\bfD^b(\A)$ which is generated by all projective modules. Show that
the composition
$$\A\lto\bfD^b(\A)\lto\bfD^b(\A)/\bfD^b(\Proj\A)$$ of canonical
functors induces an equivalence
$\bfS(\A)\to\bfD^b(\A)/\bfD^b(\Proj\A)$ of triangulated categories.

\bigskip
\noindent{\bf (7.5.1)} Let $\A$ be a Frobenius category and $\tilde\A$
the full subcategory of acyclic complexes with injective components in
$\bfC(\A)$. Show that $\tilde\A$ is a Frobenius category (with respect
to the degree-wise split exact sequences) and that the functor
$\bfS(\tilde\A)\to\bfS(\A)$ sending $X$ to $Z^0X$ is an equivalence.

\end{document}